\documentclass[preprint,12pt]{elsarticle}
\usepackage{amsfonts,amsmath,amssymb}
\usepackage{tikz}
\usepackage[colorlinks=true, citecolor=black, linkcolor=black, urlcolor=black]{hyperref}
\usepackage{todonotes}

\newcommand{\commA}[1]{\marginpar{%
\begin{color}{red}
\vskip-\baselineskip 
\raggedright\footnotesize
\itshape\hrule \smallskip A: #1\par\smallskip\hrule\end{color}}}

\newtheorem{thm}{Theorem}[section]
\newtheorem{lem}[thm]{Lemma}

\newtheorem{cor}[thm]{Corollary}
\newdefinition{df}{Definition}[section]
\newdefinition{rem}{Remark}[section]
\newdefinition{ex}{Example}[section]
\newproof{pf}{Proof}
\newproof{pot}{Proof of Theorem}
\numberwithin{equation}{section}

\journal{
}

\begin{document}

\begin{frontmatter}

\title{The height of an infinite parallelotope  is infinite
\footnote{\bf To all fearless Ukrainians defending not only their country, but the whole civilization 
against putin's 
\href{https://en.wikipedia.org/wiki/Ruscism}{rashism}
}
}
\author{A.V.~Kosyak\corref{cor1}}
\ead{kosyak02@gmail.com}
\address{Institute of Mathematics, Ukrainian National Academy of Sciences,\\
3 Tereshchenkivs'ka Str., Kyiv, 01601, Ukraine}
\address{
London Institute for Mathematical Sciences,\\
21 Albemarle St, London W1S 4BS, UK
}

\cortext[cor1]{Corresponding author}

\begin{abstract}
We show that 
$\frac{\Gamma(f_0,f_1,\dots,f_m)}
{\Gamma(f_1,\dots,f_m)}=\infty$ 
for $m+1$ vectors having the properties that no  non-trivial linear combination
of them belongs to $l_2(\mathbb N)$. This property is essential in the proof of the irreducibility of unitary representations  of some infinite-dimensional groups.
\end{abstract}

\begin{keyword}
Gram determinant \sep  parallelotope volume \sep generalized characteristic polynomial
\sep infinite-dimensional groups
 \sep irreducible representation \sep



\MSC[2020] 22E65 (5 \sep 15 \sep 26 \sep 40) 

\end{keyword}

\end{frontmatter}
\newpage
\tableofcontents

\section{One lemma about $m$ infinite vectors}
In the study of the irreducibility of unitary representations of infinite-dimensional groups  we need Lemma~\ref{l.min=proj.m}, see 
\cite{KosJFA17}--\cite{KosMor-Arx23}.
The particular case of this lemma was proved for $m=1$ in \cite{KosJFA17,Kos_B_09}
\begin{lem}
\label{l.min=proj.m}
Let $m\in \mathbb N$ and $(f_r)_{r=0}^{m}$ be $m+1$ infinite real vectors $f_r=(f_{rk})_{k\in \mathbb N},$
$0\leq r\leq m$
such that for all $\big(C_0,\dots,C_{m}\big)\in {\mathbb R}^{m+1}\setminus\{0\}$  holds
\begin{equation}
\label{norm=infty.m}
\sum_{r=0}^{m}C_rf_r\not\in l_2(\mathbb N),\quad\text{i.e.,} \quad
\sum_{k\in \mathbb N}\Big| \sum_{r=0}^{m}C_rf_{rk}\Big|^2=\infty.
\end{equation}
Denote by $f_r^{(n)}=(f_{rk})_{k=1}^n\in \mathbb R^n$ the {\rm projections} of the vectors $f_r$ on the subspace $\mathbb R^n$. Then for all 
$s$ with $0\leq s\leq m$
\begin{equation}
\label{final.m}
\frac{\Gamma(f_0,f_1,\dots,f_{m})}
{\Gamma(f_0,\dots,\hat{f_s},\dots,f_{m})}:=
\lim_{n\to\infty}\frac{\Gamma(f_0^{(n)},f_1^{(n)}\dots,f_{m}^{(n)})}{\Gamma(f_0^{(n)},\dots,\widehat{f_s^{(n)}},\dots,f_{m}^{(n)})}=\infty,
\end{equation}
where $\hat{f_s}$ means that the vector $f_s$ is absent and $\Gamma(f_0,f_1,\dots,f_{m})$ is {\rm the Gram determinant}.
\end{lem}

\begin{lem}
\label{l.det/det.m} 
Let we have $m+1$ real vectors $(f_k)_{k=0}^m$  such that
$\sum_{k=0}^mC_kf_k\not\in l_2(\mathbb N)$ for {\rm any  nontrivial combination} $(C_k)_{k=0}^m$.
Then for any $s,\,0\leq s\leq m$ 
\begin{equation}
\label{det/det.m}
\frac{{\rm det}
\big(I_{m+1}+\gamma(f_0,\dots, f_m)\big)}
{{\rm det}\big(I_{m}\!+\!\gamma(f_0,\dots,\hat{f_s},\dots,f_m)\big)}\!=\!
\lim_{n\to\infty}
\frac{
{\rm det}\big(I_{m+1}+\gamma(f_0^{(n)},\dots, f_m^{(n)})\big)}
{{\rm det}\big(I_{m}\!+\!\gamma(f_0^{(n)},\dots,\widehat{f_s^{(n)}},\dots,f_m^{(n)})\big)}
\!=\!\infty.
\end{equation}
Here $I_m\!=\!{\rm diag}(1,\dots,1)\in{\rm Mat}(m,\mathbb R,)$
and $\gamma(f_0,\dots, f_m)$ is  {\rm the Gram matrix}.
\end{lem}
\begin{pf}
The proof follows from Lemma~\ref{l.min=proj.m} and  \eqref{G.detC-LI.2}.
\qed\end{pf}
\section{The particular cases}
\subsection{Case $m=1$}
\begin{lem}[\cite{KosJFA17,Kos_B_09}]
\label{l.min=proj}
Let $f_r\!=\!(f_{rk})_{k\in{\mathbb N}},\,\,0\!\leq\! r\!\leq\! 1$ 
be two real vectors  such that
\begin{equation}
\label{norm=infty}
C_0 f_0+C_1f_1\not\in l_2(\mathbb N)\quad
\quad\text{for all}\quad
(C_0,C_1)\in {\mathbb R}^2\setminus\{0\},\quad\text{then}
\end{equation}
\begin{equation}
\label{final.}
\frac{\Gamma(f_0,f_1)}{\Gamma(f_1)}\!=\!
\lim_{n\to\infty}\frac{\Gamma(f_0^{(n)},f_1^{(n)})}{\Gamma(f_1^{(n)})}\!=\!\infty,\quad  
\frac{\Gamma(f_0,f_1)}{\Gamma(f_0)}\!=\!
\infty.
\end{equation}
\end{lem}
\begin{pf}
The initial proof can be found in  \cite{KosJFA17, Kos_B_09}.
{\it Here we give different proof that can be generalized for an arbitrary $m\in \mathbb N$}.
For $t\in \mathbb R$ and $f_0,f_1\in H$, where $H$ is some Hilbert space,  define the quadratic form
\begin{eqnarray}
\nonumber
&&F_1(t)=\Vert tf_1-f_0\Vert ^2=
t^2(f_1,f_1)-2t(f_1,f_0)+(f_0,f_0)=\\
\nonumber
&&(A_1t,t)-2(t,b)+(f_0,f_0),\,\,\,\text{where} \,\,\,b=(f_1,f_0)\in \mathbb R.\quad\text{We have}\\
\label{min.F1}
&& \min_{t\in \mathbb R}
F_1(t)=F_1(t_0)=
\frac{\Gamma(f_0,f_1)}{\Gamma(f_1)},
\quad\text{where}\quad
t_0=\frac{(f_1,f_0)}{(f_1,f_1)},
\end{eqnarray}
 and $A_1$ is the {\it Gram matrix}  (see Definition~\ref{d.Gram-det})
\begin{equation}
\label{A_1}
A_1=\gamma(f_1)=(f_1,f_1).
\end{equation}
 Consider the following matrix
\begin{equation}
\label{B(1)}
B_1:=
\gamma(f_0,f_1) 
=\left(\begin{array}{cc}
(f_0,f_0)&(f_0,f_1)\\
(f_1,f_0)&(f_1,f_1)
\end{array}\right),\quad B^{(n)}_1:=
\gamma(f^{(n)}_0,f^{(n)}_1),
\end{equation}
then 
\begin{equation}
\label{t_0.1} 
t_0=\frac{(f_1,f_0)}{(f_1,f_1)}=
\frac{-A^0_1(B_1)}{A^0_0(B_1)}.
\end{equation}
If we replace the vectors $f_0,f_1$ with 
$f^{(n)}_0,f^{(n)}_1$, the formulas \eqref{min.F1} and \eqref{t_0.1}  become
\begin{equation}
 \label{min.F10n}
 \min_{t\in \mathbb R}
F^{(n)}_{1}(t)\!=F^{(n)}_{1}(t^{(n)}_0)\!=\!
\frac{\Gamma(f^{(n)}_0,f^{(n)}_1)}{\Gamma(f^{(n)}_1)},\quad 
t^{(n)}_0\!=\!\frac{(f^{(n)}_1,f^{(n)}_0)}{(f^{(n)}_1,f^{(n)}_1)}\!=\!
\frac{-A^0_1(B^{(n)}_1)}{A^0_0(B^{(n)}_1)}.
\end{equation}
Suppose that for all $n\in \mathbb N$
\begin{equation}
\label{<C.m=1} 
\frac{\Gamma(f_0^{(n)},f_1^{(n)})}{\Gamma(f_1^{(n)})}\leq C.
\end{equation}
Without loss of generality, we can assume that  for all $n\in \mathbb N$
\begin{equation}
\label{f_0<f_1.1}
\Gamma(f^{(n)}_0)\leq C_{1}\Gamma(f^{(n)}_1),
\end{equation}
we choose some subsequence $n_k$, if necessary. Then it is sufficient to verify the first part of \eqref{final.}.  
We prove that the seguence $t^{(n)}_0$ defined by \eqref{min.F10n} is bounded. 
Since all the matrices $\gamma(f_0^{(n)},f_1^{(n)}) $ are positively defined, we have
\begin{equation*}
1\geq \frac{(f_0^{(n)},f_1^{(n)})^2}{(f_0^{(n)},f_0^{(n)})(f_1^{(n)},f_1^{(n)})} 
\stackrel{\eqref{f_0<f_1.1}}{\geq}
\frac{(f_0^{(n)},f_1^{(n)})^2}{C_{1}(f_1^{(n)},f_1^{(n)})^2}=
\frac{1}{C_{1}}\big(t^{(n)}_0\big)^2.
\end{equation*}
Hence, the sequence $t^{(n)}_0$ is bounded. Therefore,  there exists a subsequence $(t^{(n_k)}_0)_{k\in \mathbb N}$ that converges to some $t\in \mathbb R$. This contradics \eqref{<C.m=1}. Indeed
\begin{equation*}
\hskip 2 cm\lim_{n\to \infty}F^{(n)}_{1}(t)=\infty,\quad  F^{(n)}_{1}(t_0^{(n_k)})\leq C,\quad
\lim_{k\to \infty}t_0^{(n_k)}=t.\hskip 2 cm\Box
\end{equation*}
\end{pf}
\subsection{Case $m=2$}

\begin{lem}
\label{l.min=proj.3}
Let $f_0,f_1,f_2$ be three infinite real vectors $f_r=(f_{rk})_{k\in \mathbb N},$\\
$0\leq r\leq 2$
such that for all $(C_0,C_1,C_2)\in {\mathbb R}^{3}\setminus\{0\}$  holds $\sum_{r=0}^2C_rf_r\not\in l_2(\mathbb N)$.
Then for all  $r,s$ with  $0\leq r<s\leq 2$ holds
\begin{equation}
\label{final.3}
\frac{\Gamma(f_0,f_1,f_2)}{\Gamma(f_r,f_s)}:=
\lim_{n\to\infty}\frac{\Gamma(f^{(n)}_0,f^{(n)}_1,f^{(n)}_2)}{\Gamma(f^{(n)}_r,f^{(n)}_s)}
=\infty. 
\end{equation}
\end{lem}
\begin{pf}

Suppose that for all $n\in \mathbb N$
\begin{equation}
\label{<C.m=2} 
\frac{\Gamma(f_0^{(n)},f_1^{(n)},f_2^{(n)})}{\Gamma(f_1^{(n)},f_2^{(n)})}\leq C.
\end{equation}
Without loss of generality, we can assume that  for all $n\in \mathbb N$
\begin{equation}
\label{(r,s)<(1,2).2}
 \Gamma(f^{(n)}_0,f^{(n)}_1)
\leq C_{2}\Gamma(f^{(n)}_1,f^{(n)}_2),\quad 
\Gamma(f^{(n)}_0,f^{(n)}_2)
\leq 
 C_{1}\Gamma(f^{(n)}_1,f^{(n)}_2),
\end{equation}
we choose some subsequence $n_k$, if necessary. 
Then it is sufficient to verify \eqref{final.3} for $(r,s)=(1,2)$.  
Define the quadratic form $F_{2}(t)$ as follows:
\begin{eqnarray}
\nonumber
&&F_{2}(t)=\Vert \sum_{r=1}^2t_rf_r-f_0\Vert ^2=
\sum_{k,r=1}^2t_kt_r(f_k,f_r)-2\sum_{k=1}^2t_k(f_k,f_0)+(f_0,f_0)\\
\nonumber
&&=(A_2t,t)-2(t,b)+(f_0,f_0),
\end{eqnarray}
where $b=(f_k,f_0)_{k=1}^2\in \mathbb R^2$ and $A_2$ is the {\it Gram matrix}  
\begin{equation}
\label{A_2}
A_2=\gamma(f_1,f_2)=\big((f_k,f_r)\big)_{k,r=1}^2.
\end{equation}
For $t_0$ defined by $A_2t_0=b$ we have by Lemma~\ref{l.d(f,v_n)}
\begin{equation*}
{\color{blue}
F_{2}(t)}
=(A_2t,t)-2(t,b)+(f_0,f_0)
=(A_2(t-t_0),(t-t_0))+
\frac{\Gamma(f_0,f_1,f_2)}{\Gamma(f_1,f_2)},
\end{equation*}  
and therefore,
\begin{equation}
 \label{F_2(t)=(At,t).1}
{\color{blue}
F_{2}(t_0)}=
\frac{\Gamma(f_0,f_1,f_2)}{\Gamma(f_1,f_2)}.
\end{equation}
Consider the following matrix
\begin{equation}
\label{B(2)}
B_2:=
\gamma(f_0,f_1,f_2) 
=\left(\begin{array}{ccc}
(f_0,f_0)&(f_0,f_1)&(f_0,f_2)\\
(f_1,f_0)&(f_1,f_1)&(f_1,f_2)\\
(f_2,f_0)&(f_2,f_1)&(f_2,f_2)
\end{array}\right).
\end{equation}
By Cramer's rule (see Lemma~\ref{l.t_0=}) we have
\begin{equation}
\label{t_0.2} 
t_0=\frac{1}{A^0_0(B_2)}
\left(\begin{array}{c}
-A^0_1(B_2)\\
-A^0_2(B_2)
\end{array}\right).
\end{equation}
Define the quadratic form $F^{(n)}_{2}(t)$ fot $t\in \mathbb R^2$ as follows: 
\begin{eqnarray*}
\nonumber
&&F^{(n)}_{2}(t)=\Vert \sum_{r=1}^2t_rf^{(n)}_r-f^{(n)}_0\Vert ^2=
(A^{(n)}_2t,t)-2(t,b)+(f^{(n)}_0,f^{(n)}_0),
\end{eqnarray*}
where $b=(f^{(n)}_k,f^{(n)}_0)_{k=1}^2\in \mathbb R^2$ and $A^{(n)}_2$ is the {\it Gram matrix}  
\begin{equation}
\label{A_2(n)}
A^{(n)}_2=\gamma(f^{(n)}_1,f^{(n)}_2)=\big((f^{(n)}_k,f^{(n)}_r)\big)_{k,r=1}^2.
\end{equation}
By 
\eqref{F_m(t)=(At,t).1} 
we have
\begin{equation*}
F^{(n)}_{2}(t^{(n)}_0)  =
  \min_{t\in \mathbb R^2}
F^{(n)}_{2}(t)=
\frac{\Gamma(f_0^{(n)},f_1^{(n)},f_{2}^{(n)})}{\Gamma(f_1^{(n)},f_{2}^{(n)})}. 
\end{equation*}
We prove that the seguence $t^{(n)}_0$ is bounded. 
If we replace the vectors $f_0,f_1,f_2$ with $f^{(n)}_0,f^{(n)}_1,f^{(n)}_2$, we will get the following expressions:
\begin{equation}
\label{t^(n).2}
t^{(n)}_0= \frac{1}{A^0_0(B^{(n)}_2)}
\left(\begin{array}{c}
-A^0_1(B^{(n)}_2)\\
-A^0_2(B^{(n)}_2)
\end{array}\right),
\end{equation}
where 
$B^{(n)}(2)$ is defined by
\begin{equation}
\label{B^n(2)}
 B^{(n)}_2:=\!\gamma(f^{(n)}_0,f^{(n)}_1,f^{(n)}_2) \!
=\!\left(\!\!\begin{array}{ccc}
(f^{(n)}_0,f^{(n)}_0)&(f^{(n)}_0,f^{(n)}_1)&(f^{(n)}_0,f^{(n)}_2)\\
(f^{(n)}_1,f^{(n)}_0)&(f^{(n)}_1,f^{(n)}_1)&(f^{(n)}_1,f^{(n)}_2)\\
(f^{(n)}_2,f^{(n)}_0)&(f^{(n)}_2,f^{(n)}_1)&(f^{(n)}_2,f^{(n)}_2)
\end{array}\!\!\right).
\end{equation}
Since all the matrices  $B^{(n)}_2$ defined by \eqref{B^n(2)}
are positively defined,
the inverse matrices $\Big(B^{(n)}_2\Big)^{-1}$ 
are also positevely defined. We have the following expression for them
\begin{equation}
\label{B(2)^(-1)}
 \Big(B^{(n)}_2\Big)^{-1}=
\frac{1}{{\rm det}B^{(n)}_2}
 \left(\begin{array}{ccc}
A^0_0(B^{(n)}_2)&A^0_1(B^{(n)}_2)&A^0_2(B^{(n)}_2)\\
A^1_0(B^{(n)}_2)&A^1_1(B^{(n)}_2)&A^1_2(B^{(n)}_2)\\
A^2_0(B^{(n)}_2)&A^2_1(B^{(n)}_2)&A^2_2(B^{(n)}_2)
\end{array}\right). 
\end{equation}
%
%
We prove that the sequence $t^{(n)}_0$ defined by \eqref{t^(n).2}
\begin{equation}
\label{t_0.2n} 
t^{(n)}_0=
\frac{1}{A^0_0(B^{(n)}_2)}
\left(\begin{array}{c}
-A^0_1(B^{(n)}_2)\\
-A^0_2(B^{(n)}_2)
\end{array}\right)
\end{equation}
is bounded when \eqref{(r,s)<(1,2).2} holds. Set
\begin{equation}
\label{t(r,s)(n).2}
t^{(n)}_{rr}=A^r_r(B^{(n)}_2),\,\,0\leq r\leq 2,\,\,
t^{(n)}_{rs}=-A^r_s(B^{(n)}_2),\,\,\,0\leq r\not=s\leq 2.
\end{equation}
Then  $t^{(n)}_0=
\Big(\frac{t^{(n)}_{01}}{t^{(n)}_{00}},\frac{t^{(n)}_{02}}{t^{(n)}_{00}}
\Big)$.
Since the matrix $\Big(B^{(n)}_2\Big)^{-1}$
is positively defined and
 \eqref{(r,s)<(1,2).2} holds,  we have
\begin{eqnarray}
\label{A:>}
&&
\big(t^{(n)}_{01}\big)^2\leq t^{(n)}_{00}t^{(n)}_{11},\quad \big(t^{(n)}_{02}\big)^2\leq t^{(n)}_{00}t^{(n)}_{22},\\
&&
\label{A:(r,s)<(1,2).2}
t^{(n)}_{22}
\leq C_2t^{(n)}_{00},\quad
t^{(n)}_{11}
\leq
C_1t^{(n)}_{00},\\
&&
\Vert t^{(n)}_0\Vert^2:=
\frac{\vert t^{(n)}_{01} \vert^2+\vert t^{(n)}_{02} \vert^2}
{\vert t^{(n)}_{00}\vert^2}\stackrel{\eqref{A:>}}{\leq}
\frac{ t^{(n)}_{11}+ t^{(n)}_{22}}{ t^{(n)}_{00}}
\\
&&
\stackrel{\eqref{A:(r,s)<(1,2).2}}{\leq}
\frac{ C_1t^{(n)}_{00}+C_2 t^{(n)}_{00}}{ t^{(n)}_{00}}=C_1+C_2.
\end{eqnarray}
Hence, the sequence $t^{(n)}_0\in \mathbb R^2$ is bounded. Therefore,  there exists a subsequence $(t^{(n_k)}_0)_{k\in \mathbb N}$ that converges to some $t\in \mathbb R^2$. This contradics \eqref{<C.m=2}. Indeed
\begin{equation*}
\hskip 1.5 cm 
\lim_{n\to \infty}F^{(n)}_{2}(t)=\infty,\quad  F^{(n)}_{2}(t_0^{(n_k)})\leq C,\quad
\lim_{k\to \infty}t_0^{(n_k)}=t.
\hskip 1.5 cm \Box
\end{equation*}
\end{pf}

%
\section{The proof of  Lemma~\ref{l.min=proj.m}}

\begin{pf} [Proof of the Lemma~\ref{l.min=proj.m}]

Suppose that for all $n\in \mathbb N$ we have
\begin{equation}
\label{F_n.m<C} 
\frac{\Gamma(f_0^{(n)},f_1^{(n)}\dots,f_{m}^{(n)})}{\Gamma(f_1^{(n)},f_2^{(n)}\dots,f_{m}^{(n)})}\leq C.
\end{equation}
Without loss of generality, we can assume that  for all $n\in \mathbb N$ and $0\leq s\leq m-1$ hold
\begin{equation}
\label{Gam<Gam.m}
\Gamma(f_0^{(n)},\dots,\widehat{f_{s}^{(n)}},\dots,f_{m}^{(n)})<C_s\Gamma(f_1^{(n)},f_2^{(n)},\dots,f_{m}^{(n)}).
\end{equation}
we choose some subsequence $n_k$, if necessary. 
Consider the following quadratic forms for $t=(t_r)_{r=1}^m\in \mathbb R^m$
\begin{equation}
\label{F_n.m}
F^{(n)}_{m}(t)
=\Vert f_0^{(n)}-\sum_{r=1}^mt_rf_r^{(n)}\Vert^2,
 \end{equation}
The forms $ F^{(n)}_{m}(t)$ 
defined by \eqref{F_n.m} 
have the following properties, for any fixed $t\in \mathbb R^m$ we have:\\
1)  $ F^{(n)}_{m}(t)\leq  F^{(n+1)}_{m}(t)$,\\
2)  $\lim_{n\to \infty}  F^{(n)}_{m}(t)=\infty$.\\
By Lemma~\ref{l.d(f,v_n)} there exists some 
$t^{(n)}_0
\in
\mathbb R^m$ such that 
\begin{equation}
\label{min-Q_n.m}
 F^{(n)}_{m,0}(t^{(n)}_0)=
  \min_{t\in \mathbb R^m}
 F^{(n)}_{m,0}(t)=
\frac{\Gamma(f_0^{(n)},f_1^{(n)},\dots,f_{m}^{(n)})}{\Gamma(f_1^{(n)},f_{2}^{(n)}\dots,f_{m}^{(n)})}. 
\end{equation}
We prove that the seguence $t^{(n)}_0$ is bounded.
To find $t^{(n)}_0$  explicitely in 
\eqref{min-Q_n.m} we 
introduce some notations.
For $t\in \mathbb R^m$ and $f_0,f_1,\dots,f_m\in H$ we define the function
\begin{eqnarray*}
F_m(t)&=&\Vert\sum_{k=1}^mt_kf_k-f_0\Vert ^2=
\sum_{k,r=1}^mt_kt_r(f_k,f_r)-2\sum_{k=1}^2t_k(f_k,f_0)+(f_0,f_0)\\
\nonumber
&=&(A_mt,t)+2(t,b)+(f_0,f_0),
\end{eqnarray*}
where $b=(f_k,f_0)_{k=1}^m\in \mathbb R^m$ and $A_m$ is the 
{\it Gram matrix}, see Definition~\ref{d.Gram-det}: 
\begin{equation}
\label{A_m}
A_m=\gamma(f_1,\dots, f_m)=\big((f_k,f_r)\big)_{k,r=1}^m.
\end{equation}
The minimum of $F_m(t)$ is attained at $t_0$ defined by  $A_mt_0=b$. By  Remark~\ref{r.At=b}, 
\eqref{At=b}  and \eqref{d^2} we get (for details see  Section~
\ref{s.3.1})
\begin{equation*}
F_m(t)
=(A_mt,t)-2(t,b)+(f_0,f_0)
=(A_m(t-t_0),(t-t_0))+
\frac{\Gamma(f_0,f_1,\dots,f_m)}{\Gamma(f_1,\dots,f_m)},
\end{equation*}  
\begin{equation}
 \label{F_m(t)=(At,t).1}
\text{and therefore,}\quad F_m(t_0)=
\frac{\Gamma(f_0,f_1,\dots,f_m)}{\Gamma(f_1,\dots,f_m)}.
\end{equation}
Consider the following matrix
\begin{equation}
\label{B(m)}
B_m:=
\gamma(f_0,f_1,\dots,f_m) 
=\left(\begin{array}{ccccc}
(f_0,f_0)&(f_0,f_1)&(f_0,f_2)&\dots&(f_0,f_m)\\
(f_1,f_0)&(f_1,f_1)&(f_1,f_2)&\dots&(f_1,f_m)\\
(f_2,f_0)&(f_2,f_1)&(f_2,f_2)&\dots&(f_2,f_m)\\
&&&\dots&\\
(f_m,f_0)&(f_m,f_1)&(f_m,f_2)&\dots&(f_m,f_m)
\end{array}\right).
\end{equation}
By Cramer's rule (see Lemma~\ref{l.t_0=}) the solution of $A_mt_0=b$ is as follows:
\begin{equation}
\label{t_0.m} 
t_0\!=\!\frac{1}{A^0_0(B_m)}
\left(\!\begin{array}{c}
-A^0_1(B_m)\\
-A^0_2(B_m)\\
\dots\\
-A^0_m(B_m)
\end{array}\!\right).
\end{equation}
If we replace the vectors $f_0,f_1,\dots,f_m$ with $f^{(n)}_0,f^{(n)}_1,\dots,f^{(n)}_m$, we will get the following expression
\begin{equation}
\label{t^(n).m}
t^{(n)}_0\!=\!
\frac{1}{A^0_0(B^{(n)}_m)}
\left(\!\!\begin{array}{c}
-A^0_1(B^{(n)}_m)\\
\dots\\
-A^0_m(B^{(n)}_m)
\end{array}\!\!\right),
\end{equation}
wher $B^{(n)}(m)$ is defined by
\begin{equation}
\label{B^n(m)}
 B^{(n)}_m:=\!\gamma(f^{(n)}_0,\dots,f^{(n)}_m) \!
=\!\left(\!\!\begin{array}{cccc}
(f^{(n)}_0,f^{(n)}_0)&(f^{(n)}_0,f^{(n)}_1)&\dots&(f^{(n)}_0,f^{(n)}_m)\\
(f^{(n)}_1,f^{(n)}_0)&(f^{(n)}_1,f^{(n)}_1)&\dots&(f^{(n)}_1,f^{(n)}_m)\\
&&\dots&\\
(f^{(n)}_m,f^{(n)}_0)&(f^{(n)}_m,f^{(n)}_1)&\dots&(f^{(n)}_m,f^{(n)}_m)
\end{array}\!\!\right).
\end{equation}
Since all the matrices  $B^{(n)}_m$ defined by \eqref{B^n(m)}
are positively defined,
the inverse matrices $\Big(B^{(n)}_m\Big)^{-1}$ 
are also positevely defined. We have the following expression for them
\begin{equation}
\label{B(m)^(-1)}
 \Big(B^{(n)}_m\Big)^{-1}\!\!\!=\!
\frac{1}{{\rm det}B^{(n)}_m}
 \left(\!\!\begin{array}{cccc}
A^0_0(B^{(n)}_m)&A^0_1(B^{(n)}_m)&\dots&A^0_m(B^{(n)}_m)\\
A^1_0(B^{(n)}_m)&A^1_1(B^{(n)}_m)&\dots&A^1_m(B^{(n)}_m)\\
&&\dots&\\
A^m_0(B^{(n)}_m)&A^m_1(B^{(n)}_m)&\dots&A^m_m(B^{(n)}_m)
\end{array}\!\!\right)\!. 
\end{equation}
We prove that the sequence 
\begin{equation}
\label{t_0.2m} 
t^{(n)}_0=
\frac{1}{A^0_0(B^{(n)}_m)}
\left(\begin{array}{c}
-A^0_1(B^{(n)}_m)\\
\dots\\
-A^0_m(B^{(n)}_m)
\end{array}\right)
\end{equation}
is bounded when \eqref{Gam<Gam.m} 
holds. Set
\begin{equation}
\label{t(r,s)(n).m}
t^{(n)}_{rr}=A^r_r(B^{(n)}_m),\,\,0\leq r\leq m,\,\,
t^{(n)}_{rs}=-A^r_s(B^{(n)}_m),\,\,\,0\leq r\not=s\leq m.
\end{equation}
Then  $t^{(n)}_0=\Big(\frac{t^{(n)}_{01}}{t^{(n)}_{00}},\frac{t^{(n)}_{02}}{t^{(n)}_{00}},\dots
\frac{t^{(n)}_{0m}}{t^{(n)}_{00}}
\Big)$.
Since the matrix $\Big(B^{(n)}_m\Big)^{-1}$
is positively defined and
 \eqref{(r,s)<(1,2).2} holds,  we have 
\begin{eqnarray}
\label{A:>.m}
&&
\big(t^{(n)}_{rs}\big)^2\leq t^{(n)}_{rr}t^{(n)}_{ss},\quad\text{for all}\quad 0\leq r<s\leq m,\\
&&
\label{A:(r,s)<(1,2).m}
t^{(n)}_{ss}\leq C_st^{(n)}_{00},\quad
\text{for all}\quad 0\leq s\leq m-1,\\
&&
\Vert t^{(n)}_0\Vert^2:=
\frac{\sum_{s=1}^{m}\vert t^{(n)}_{0s} \vert^2}
{\vert t^{(n)}_{00}\vert^2}\stackrel{\eqref{A:>.m}}{\leq}
\frac{\sum_{s=1}^{m} t^{(n)}_{ss}}{ t^{(n)}_{00}}
\\
&&
\stackrel{\eqref{A:(r,s)<(1,2).m}}{\leq}
\frac{\sum_{s=1}^{m}C_st^{(n)}_{00}
}{ t^{(n)}_{00}}=\sum_{s=1}^{m}C_s.
\end{eqnarray}
Hence, the sequence $t^{(n)}_0\in \mathbb R^m$ is bounded. Therefore,
 there exists a subsequence $(t^{(n_k)}_0)_{k\in \mathbb N}$ that converges to some $t\in \mathbb R^m$. This contradics \eqref{F_n.m<C}. Indeed
\begin{equation*}
\hskip 1.5 cm 
\lim_{n\to \infty}F^{(n)}_{m,0}(t)=\infty,\quad  F^{(n)}_{m,0}(t_0^{(n_k)})\leq C,\quad
\lim_{k\to \infty}t_0^{(n_k)}=t.\hskip 1.5 cm \Box
\end{equation*}
\end{pf}

\section{How far is a vector from a hyperplane?}

\subsection{The distance of a vector from a hyperplane}
\label{s.3.1}
In this section we follow  \cite[Section 1.3]{Kos-hpl-arx23}. We start with a classical result, see, e.g.
\cite{Gan58}.
Consider the hyperplan $V_n$ generated by $n$ arbitrary vectors $f_1,\dots, f_{n}$ in some Hilbert space  $H$. 
\begin{lem} 
\label{l.d(f,v_n)}
The square of the distance $d(f_0,V_n)$ of a vector $f_0$ from the hyperplaen $V_n$ is given by the ratio of two Gram determinants (see Definition \ref{d.Gram-det})
\begin{equation}
\label{d(f,v_n)}
d^2(f_0,V_n)=
\frac{\Gamma(f_0,f_1,f_2,\dots, f_n,)}{\Gamma(f_1,f_2,\dots, f_n)}.
\end{equation}
\end{lem}
\begin{pf} We follow closely the book by Axiezer and Glazman \cite{AhiGlaz93}.
Set $f=\sum_{k=1}^nt_kf_k\in V_n$ and $h=f-f_0$. Since $h$ should be orthogonal to $V_n$ we conclude that
$f_r\perp h$,  i.e., $(f_r,h)=0$ for all $r$, or
\begin{equation}
 \label{f_r-perp-h}
\sum_{k=1}^nt_k(f_r,f_k)=(f_r,f_0),\quad 1\leq r\leq n. 
\end{equation}
Set $A=\gamma(f_1,f_2,\dots,f_n)$ and $b=(f_k,f_0)_{k=1}^n\in \mathbb R^{n} $. By definition we have 
\begin{equation}
\label{d^2} 
d^2=\min_{f\in V_n}\Vert f-f_0\Vert^2=(At,t)-2(t,b)+(f_0,f_0).
\end{equation}
Since $d^2=(h,h)=(f_0,h)$ we conclude that $d^2=\sum_{k=1}^nt_k(f_0,f_k)-(f_0,f_0)$ or
\begin{equation}
 \label{d^2=(f_0,h)}
\sum_{k=1}^nt_k(f_0,f_k)=(f_0,f_0)-d^2. 
\end{equation}
So we have the system of equations:
\begin{equation}
\label{syst}
\left\{
\begin{array}{ccc}
t_1(f_1,f_1)+t_2(f_1,f_2)+\dots +t_n(f_1,f_n)&=&(f_1,f_0)\\
t_1(f_2,f_1)+t_2(f_2,f_2)+\dots +t_n(f_2,f_n)&=&(f_2,f_0)\\
\dots&&\\
t_1(f_n,f_1)+t_2(f_n,f_2)+\dots +t_n(f_n,f_n)&=&(f_n,f_0)\\
t_1(f_0,f_1)+t_2(f_0,f_2)+\dots +t_n(f_0,f_n)&=&(f_0,f_0)-d^2
\end{array}
\right..
\end{equation}
Excluding $t_k$ from the system we get $d^2=\frac{\Gamma(f_0,f_1,f_2,\dots, f_n,)}{\Gamma(f_1,f_2,\dots, f_n)}$.
\qed\end{pf}
\begin{rem}
\label{r.At=b}
From the system \eqref{syst} we conclude that $At=b$, where $b=(f_k,f_0)_{k=1}^n\in \mathbb R^{n} $, hence $t=A^{-1}b$. 
By \eqref{d^2} we get 
\begin{equation}
\label{At=b} 
d^2=(f_0,f_0)-(A^{-1}b,b)=\frac{\Gamma(f_0,f_1,f_2,\dots, f_n)}{\Gamma(f_1,f_2,\dots, f_n)}.
\end{equation}
See also \cite[Chap. 4.3, Lemma 4.3.2]{Kos_B_09}.
\end{rem}
\subsection{Cramer's rule reformulated
}

Consider two 
{\it Gram matrices}, 
(see Definition~\ref{d.Gram-det}): 
\begin{equation}
\label{A_m,B_m}
A_m=\gamma(f_1,\dots, f_m)=\big((f_k,f_r)\big)_{k,r=1}^m,\quad B_m=
\gamma(f_0,f_1,\dots,f_m), 
\end{equation}
and a vector $b=(f_k,f_0)_{k=1}^m\in \mathbb R^m$. The solution of the equation $A_mt=b$ is  as follows.
\begin{lem}
\label{l.t_0=}
We have
\begin{equation}
 \label{t-lin-eq-(A)}
 t=A_m^{-1}b=\frac{1}{A^0_0(B_m)}
\left(\!\begin{array}{c}
-A^0_1(B_m)\\
-A^0_2(B_m)\\
\dots\\
-A^0_m(B_m)
\end{array}\!\right).
\end{equation}
\end{lem}
\begin{pf}
By {\it Cramer's rule} if we have a system of linear equations
\begin{equation}
\label{lin-eq}
At=b
\end{equation}
where $A\in {\rm Mat}(m,\mathbb C)$ with ${\rm det }\,A\not=0$ and $t,b\in \mathbb C^m$, than the soultions are given by the following formulas:
\begin{equation}
\label{t-lin-eq}
t_k=\frac{{\rm det}\,(A_k)}{{\rm det}\,(A)},\quad 1\leq k\leq m,
\end{equation}
where 
$A_{k}$ is the matrix formed by replacing the $i$-th column of $A$ by the column vector $b$.
Consider the  matrix $B_m$ defined by \eqref{A_m,B_m}.
We have 
\begin{equation}
 \label{A(r,0)=A(r)}
 {\rm det}\,(A)=A^0_0(B_m),\quad{\rm det}\,(A_k)=
 -A^0_k(B_m),\quad 1\leq k\leq m,
\end{equation}
thus implying \eqref{t-lin-eq-(A)}.
Recall that $A^r_s(C)$ denote {\it cofactors} of the matrix $C$, for details, see Section~\ref{sec.gen.har.pol}.
\qed\end{pf}

\subsection{Some estimates}
\label{s.1.4.1}
%
%
We use some material from \cite{Kos_B_09}, Section 1.4.1, pp. 24-25.
\begin{lem} 
\label{1.l.min}
 For a strictly positive operator $A$ (i.e., $(Af,f)>0,\,\,f\not=0$) acting in ${\mathbb R}^n$ and a vector
$b\in{\mathbb R}^n\backslash\{0\}$ we have
\begin{equation}
\label{A.min3} \min_{x\in{\mathbb R}^n}\Big((Ax,x)\mid
(x,b)=1\Big)=(A^{-1}b,b)^{-1}.
\end{equation}
The minimum is assumed  for
$x=A^{-1}b\left((A^{-1}b,b)\right)^{-1}$.
\end{lem}
Lemma \ref{1.l.min} is a direct generalization of the well known
result (see, for example, \cite{BecBel61}, Chap. I, \S 52), for $a_k>0,\,\,1\leq k \leq n$ we have
\begin{equation}
\label{A.min1}
 \min_{x\in{\mathbb
R}^n}\Big(\sum_{k=1}^na_kx_k^2\mid \sum_{k=1}^n x_k=1\Big)=
\Big(\sum_{k=1}^n \frac{1}{a_k}\Big)^{-1}.
\end{equation}
We will also use  the same result   in a slightly different
form:
\begin{equation}
\label{A.min2}
 \min_{x\in{\mathbb
R}^n}\Big(\sum_{k=1}^na_kx_k^2\mid \sum_{k=1}^n x_kb_k=1\Big)=
\Big(\sum_{k=1}^n \frac{b_k^2}{a_k}\Big)^{-1}.
\end{equation}
The minimum is  assumed   for
$x_k=\frac{b_k}{a_k}\Big(\sum_{k=1}^n\frac{b_k^2}{a_k}\Big)^{-1}$.

%

\subsection{Gram determinants and Gram matrices}
\begin{df}
\label{d.Gram-det}
Let us recall the definition of
	the Gram  determinant and the Gram  matrix   (see
	\cite{Gan58}, Chap IX, \S 5). Given the  vectors $x_1,x_2,..., x_m$ in some
	Hilbert space $H$ the {\it Gram
		matrix} $\gamma(x_1,x_2,..., x_m)$ is defined  by the formula
$$
\gamma(x_1,x_2,..., x_m)=\big((x_k,x_n)\big)_{k,n=1}^m.
	$$
	\index{Gram determinant}\index{Gram matrix}
	The determinant of this matrix is called the {\it Gram  determinant} for
	the vectors $x_1,x_2,..., x_m$ and is denoted by
	$\Gamma(x_1,x_2,..., x_m)$.
	 Thus,
	 %
 %
	\begin{equation}
	\label{Gram-det}
	\Gamma(x_1,x_2,\dots, x_m):={\rm det}\,\gamma(x_1,x_2,\dots, x_m).
	\end{equation}
\end{df}
\subsection{The generalized characteristic polynomial
and its properties}
\label{sec.gen.har.pol}
\index{polynomial!characteristic!generalized}
{\bf Notations}. For a matrix $C\in {\rm Mat}(n,\mathbb R)$  and $1\leq i_1<i_2<\dots i_r\leq n,$ 
$1\leq j_1<j_2<\dots j_r\leq n,\,\,r\leq n$
denote by 
$$
M^{i_1i_2\dots i_r}_{j_1j_2\dots j_r}(C)\quad\text{ and} \quad
A^{i_1i_2\dots i_r}_{j_1j_2\dots j_r}(C)
$$
the corresponding {\it minors} and {\it cofactors} of the matrix $C$.
\begin{df}{\rm (\cite[Ch.1.4.3]
{Kos_B_09})}
\label{d.G_k(lambda)} 
For the matrix $C\in {\rm
Mat}(m,{\mathbb C})$ and $\lambda =(\lambda_1,\dots,\lambda_m)\in {\mathbb C}^m$ define  the {\it
generalization of the characteristic polynomial}, $p_C(t)={\rm
det}\,(tI-C),\, t\in {\mathbb C}$  
as follows:
\begin{equation}
\label{P_C(lambda)}
P_C(\lambda)={\rm det}\,C(\lambda),\quad\text{where}\quad
C(\lambda)=
{\rm diag}(\lambda_1,\dots,\lambda_m)
+C.
\end{equation}
\end{df}
\begin{lem}
\label{l.detC-LI}
{\rm (\cite[Ch.1.4.3]{Kos_B_09})}
For the generalized characteristic polynomial $P_C(\lambda)$
of
$C\!\in\!{\rm Mat}(m,{\mathbb C})$ and
 $\lambda=(\lambda_1,\lambda_2,...,\lambda_m)\in {\mathbb
C}^m$ we have 
\begin{equation}
\label{detC-LI}
P_C(\lambda)=
{\rm det}\,C+
\sum_{r=1}^m\sum_{1\leq i_1<i_2<...<i_r\leq
m}\lambda_{i_1}\lambda_{i_2}...\lambda_{i_r}A^{i_1i_2...i_r}_{i_1i_2...i_r}(C).
\end{equation}
\end{lem}
\begin{rem}
\label{r.P_C(lam)}
If we set
$\lambda_\alpha=\lambda_{i_1}\lambda_{i_2}...\lambda_{i_r}$, where
$\alpha=(i_1,i_2,...,i_r)$ and
$A^\alpha_\alpha(C)=A^{i_1i_2...i_r}_{i_1i_2...i_r}(C),\,\,
M^\alpha_\alpha(C)=M^{i_1i_2...i_r}_{i_1i_2...i_r}(C),\,\,
\lambda_\emptyset=1,\,\,A^\emptyset_\emptyset(C)={\rm det}\,C$ we
may write (\ref{detC-LI}) as follows:
\begin{equation}
\label{A.detC-LI.2} 
P_C(\lambda)={\rm det}\,C(\lambda)=
\sum_{\emptyset\subseteq\alpha\subseteq\{1,2,...,m\}}\lambda_\alpha
A^\alpha_\alpha(C),
\end{equation}
\begin{equation}
\label{M.detC-LI.2}
P_C(\lambda)={\rm det}\,C(\lambda)=
\left(\prod_{k=1}^n\lambda_k\right)\sum_{\emptyset\subseteq\alpha\subseteq\{1,2,...,m\}}
\frac{M^\alpha_\alpha(C)}{\lambda_\alpha},
\end{equation}
\end{rem}

Let 
\begin{equation}
\label{X(mn)}
X=X_{mn} =\left(
\begin{array}{cccc}
x_{11}&x_{12}&...&x_{1n}\\
x_{21}&x_{22}&...&x_{2n}\\
...   &...   &...&...\\
x_{m1}&x_{m2}&...&x_{mn}
\end{array}
\right).
\end{equation}
Setting 
%
%
\begin{equation}
\label{x_k,y_r=}
x_k=(x_{1k},x_{2k},...,x_{mk})\in{\mathbb R}^m,
\quad
y_r=(x_{r1},x_{r2},...,x_{rn})\in{\mathbb R}^n,
\end{equation}
we get
\begin{equation}
\label{X^*X}
X^*X=
 \left(
\begin{array}{cccc}
(x_1,x_1)&(x_1,x_2)&...&(x_1,x_n)\\
(x_2,x_1)&(x_2,x_2)&...&(x_2,x_n)\\
...   &...   &...&...\\
(x_n,x_1)&(x_n,x_2)&...&(x_n,x_n)
\end{array}
\right)=\gamma(x_1,x_2,..., x_n),
\end{equation}
\begin{equation}
\label{XX^*}
XX^*= \left(
\begin{array}{cccc}
(y_1,y_1)&(y_1,y_2)&...&(y_1,y_m)\\
(y_2,y_1)&(y_2,y_2)&...&(y_2,y_m)\\
...   &...   &...&...\\
(y_m,y_1)&(y_m,y_2)&...&(y_m,y_m)
\end{array}
\right)=\gamma(y_1,y_2,..., y_m),
\end{equation}
therefore, we obtain
\begin{equation}
 \label{d(X^*X)=d(XX^*)}\
\Gamma(x_1,x_2,..., x_n)={\rm det}(X^*X)={\rm det}(XX^*)=\Gamma(y_1,y_2,..., y_m).
\end{equation}

\section{The explicit expression for $C^{-1}(\lambda)$ and  $(C^{-1}(\lambda)a,a)$}
In this section we follow \cite[Section~2]{Kos-hpl-arx23}.
Fix $C\in {\rm Mat}(n,\mathbb R),$ $a\in \mathbb R^n$ and $\lambda\in \mathbb C^n$. 
Our aim is to find the explicit formulas for  
$C^{-1}(\lambda)$ and  $(C^{-1}(\lambda)a,a)$, where $C(\lambda)$ is defined by \eqref{P_C(lambda)}.  
Set $M(i_1i_2\dots i_r)(C)$
$=M^{i_1i_2\dots i_r}_{i_1i_2\dots i_r}(C)$ and 
$a_{i_1i_2\dots i_r}=(a_{i_1},a_{i_2},\dots,a_{i_r})$.
Let also $C_{i_1i_2\dots i_r}$ be the corresponding {\it submatrix} of the matrix $C$.
The elements of this matrix are on the intersection of $i_1,i_2,\dots, i_r$ rows and column of the matrix $C$.
Denote by $A(C_{i_1i_2\dots i_r})$ the matrix of the cofactors of the first order of the matrix $C_{i_1i_2\dots i_r}$, i.e. 
\begin{equation}
\label{A(C)}
A(C_{i_1i_2\dots i_r})=(A^{i}_j(C_{i_1i_2\dots i_r}))_{1\leq i,j\leq r} 
\end{equation}

Let $n=3$, then $A(C_{123})=A(C)$ is the following matrix: 
\begin{equation}
 \label{A(C)3}
A(C)=A(C_{123})=
\left(
\begin{array}{ccc}
A^1_1&A^1_2&A^1_3\\
A^2_1&A^2_2&A^2_3\\
A^3_1&A^3_2&A^3_3
\end{array}
\right)=
\left(
\begin{array}{ccc}
M^{23}_{23}&-M^{23}_{13}&M^{23}_{12}\\
-M^{13}_{23}&M^{13}_{13}&-M^{13}_{12}\\
M^{12}_{23}&-M^{12}_{13}&M^{12}_{12}
\end{array}
\right),
\end{equation}
%
%
where we write $M^{ij}_{rs}$ instead of $M^{ij}_{rs}(C)$ and $A^i_j$ instead of $A^i_j(C)$.
\begin{rem}
\label{r.C_{12}=}
Let $A^T$ be the transposed matrix
of $A$. Then
\begin{equation}
 \label{A(C_{12})=}
A^T(C_{i_1i_2\dots i_r})={\rm det}\,C^{-1}_{i_1i_2\dots i_r}\Big(C^{-1}_{i_1i_2\dots i_r}\Big),
\end{equation}

\end{rem}
In what follows we will  consider the submatrix $C_{i_1i_2\dots i_r}$ of the matrix $C\in {\rm Mat}(n,\mathbb R)$ as an appropriate element of ${\rm Mat}(n,\mathbb R)$.
\begin{thm}
\label{t.(C^{-1}a,a)}
 For the matrix $C(\lambda)$ defined by  \eqref{P_C(lambda)} 
 $a\in \mathbb R^n$  and $\lambda\in \mathbb C^n$ we have 
 {\small
   \begin{equation}
 \label{Del_n(lam,C)}
 P_C(\lambda)=
 \Big(\prod_{k=1}^n\lambda_k\Big)\sum_{r=1}^n
\sum_{1\leq i_1<i_2<\dots <i_r\leq n}
\frac{M(i_1i_2\dots i_r)}{\lambda_{i_1}\lambda_{i_2}\dots \lambda_{ i_r}},
\end{equation}
 \begin{equation}
 \label{C^{-1}(lam)}
C^{-1}(\lambda)=\frac{1}{P_C(\lambda)}\Big(\prod_{k=1}^n\lambda_k\Big)
\sum_{r=1}^n
\sum_{1\leq i_1<i_2<\dots i_r\leq n}
\frac{A(C_{i_1i_2\dots i_r})}{\lambda_{i_1}\lambda_{i_2}\dots \lambda_{ i_r}},
 \end{equation} 
 \begin{equation}
 \label{(C^{-1}(lam)a,a)}
(C^{-1}(\lambda)a,a)=\frac{1}{P_C(\lambda)}\Big(\prod_{k=1}^n\lambda_k\Big)
\sum_{r=1}^n
\sum_{1\leq i_1<i_2<\dots i_r\leq n}
\frac{(A(C_{i_1i_2\dots i_r})a_{i_1i_2\dots i_r},a_{i_1i_2\dots i_r})}{\lambda_{i_1}\lambda_{i_2}\dots \lambda_{ i_r}}.
 \end{equation}
 }
\end{thm}

\subsection{The case where $C$ is the Gram matrix}
Fix the matrix $X_{mn}$ defined by \eqref{X(mn)}. Denote by $C$ the Gram matrix 
$\gamma(x_1,x_2,..., x_n)$, 
i.e.,
\begin{equation}
\label{C=gamma}
C=\gamma(x_1,x_2,..., x_n),
\end{equation}
where $(x_1,x_2,..., x_n)$ are defined by \eqref{x_k,y_r=} 
and $\gamma(x_1,x_2,..., x_n)$  by \eqref{X^*X}.
In what follows we consider the operator $C(\lambda)$  defined  
by \eqref{P_C(lambda)}.
\begin{rem}
 \label{r.P_C(lam).7}
In this case  we have
%
%
\begin{eqnarray}
 \label{G.detC-LI.2} 
P_C(\lambda)&=&
{\rm det}\Big(\sum_{k=1}^n\lambda_kE_{kk}+\gamma(x_1,x_2,...,
x_n)\Big)\\
\nonumber
&=&\prod_{k=1}^n\lambda_k\Big(1+\sum_{r=1}^n\sum_{1\leq
	i_1<i_2<...<i_r\leq
	m}\Big(\lambda_{i_1}\lambda_{i_2}...\lambda_{i_r}\Big)^{-1}
\Gamma(x_{i_1},x_{i_2},...,x_{i_r})\Big)\\
\nonumber
&=&\prod_{k=1}^n\lambda_k\Big(1+\sum_{r=1}^{n} \sum_{
\substack{1\leq i_1<i_2<...<i_r\leq n;\\
1\leq j_1<j_2<...<j_r\leq n}
}
\Big(\lambda_{i_1}\lambda_{i_2}...\lambda_{i_r}\Big)^{-1}
\Big(M^{i_1i_2...i_r}_{j_1j_2...j_r}(X)\Big)^2\Big),
\end{eqnarray}
where we have used the following formula (see \cite{Gan58}, Chap IX, \S
5 formula (25)):
\begin{equation}
\label{Gramm(x,y)=M^2(X)}
\Gamma(x_{i_1},x_{i_2},...,x_{i_r})= \sum_{1\leq
	j_1<j_2<...<j_r\leq m}
\left(M^{i_1i_2...i_r}_{j_1j_2...j_r}(X)\right)^2.
\end{equation}
\end{rem}

\subsection{Case $m=2$}

Fix two natural numbers $n,m\in \mathbb N$ with $m\leq n$, two matrices $A_{mn}$ and $X_{mn}$, vectors $g_k\in\mathbb R^{m-1},\,\,1\leq k\leq n$
and $a\in \mathbb R^n$ as follows
\begin{equation} 
\label{A(mn)}
A_{mn}\!=\!\left(
	\begin{array}{cccc}
	a_{11}&a_{12}          &...&a_{1n}\\
	a_{21}&a_{22}          &...&a_{2n}\\
	&&...&\\
	a_{m1}&a_{m2}          &...&a_{mn}
	\end{array}
	\right),\,\, 
g_k\!=\!\left(
	\begin{array}{c}
	a_{2k}\\
	a_{3k}\\
	... \\
	a_{mk}
	\end{array}
	\right)\in \mathbb R^{m-1},\,\,\, a\!=\!(a_{1k})_{k=1}^n\in \mathbb R^n.	
\end{equation}
Set $C=\gamma(g_1,g_2,\dots,g_n)$. We calculate $\Delta_n(\lambda, C)$ and $(C^{-1}(\lambda)a,a)$ for an arbitrary $n$. 
Consider the matrix \eqref{X(mn)}
\begin{equation}
\label{X(mn).1}
X_{mn} \!=\!\left(
\begin{array}{cccc}
x_{11}&x_{12}&...&x_{1n}\\
x_{21}&x_{22}&...&x_{2n}\\
&&...&\\
x_{m1}&x_{m2}&...&x_{mn}
\end{array}
\right),\quad \text{where}\quad x_{kr}\!=\!\frac{a_{1k}}{\sqrt \lambda_k},\quad y_r\!=\!(x_{rk})_{k=1}^n\in \mathbb R^n.
\end{equation}
%
%
\begin{lem}[\cite{Kos-hpl-arx23}, Lemma~2.2]
\label{l.m=2}
For $m=2$ we have
\begin{equation}
\label{m=2}
(C^{-1}(\lambda)a,a)=\Delta(y_1,y_2):=\frac{\Gamma(y_1)+\Gamma(y_1,y_2)}{1+\Gamma(y_2)},
\end{equation}
where $y_1$ and $y_2$ are defined as follows
\begin{equation}
\label{y(1),y(2)=} 
y_1=\left(\frac{a_{1k}}{\sqrt{\lambda_k}}\right)_{k=1}^n\quad y_2=\left(\frac{a_{2k}}{\sqrt{\lambda_k}}\right)_{k=1}^n.
\end{equation} 
\end{lem}

\subsection{Case $m=3$
}
Define
$$
\Delta(y_1,y_2,y_3):=\frac{\Gamma(y_1)+\Gamma(y_1,y_2)+\Gamma(y_1,y_3)+\Gamma(y_1,y_2,y_3)}{1+\Gamma(y_2)+\Gamma(y_3)+\Gamma(y_2,y_3)}.
$$
\begin{lem}[\cite{Kos-hpl-arx23}, Lemma~2.3]
\label{l.m=3}
For $m=3$ we have
\begin{equation}
\label{m=3}
(C^{-1}(\lambda)a,a)=\Delta(y_1,y_2,y_3),
\end{equation}
where the $y_r$ are defined as follows:
\begin{equation}
\label{f,g,h=} 
y_r=
\left(\frac{a_{rk}}{\sqrt{\lambda_k}}\right)_{k=1}^n\in \mathbb R^n,\quad 1\leq r\leq 3.
\end{equation} 
\end{lem}

\subsection{General $m$}
 We note by \eqref{G.detC-LI.2} that
\begin{eqnarray*}
&& 
\Delta(y_1,y_2)\!=\!\frac{{\rm det}
\big(I_2\!+\!\gamma(y_1,y_2)\big)}
{{\rm det}\big(I_{1}+\gamma(y_2)\big)}\!-\!1,\,\,
\Delta(y_1,y_2,y_3)\!=\!\frac{{\rm det}
\big(I_3\!+\!\gamma(y_1,y_2,y_3)\big)}
{{\rm det}\big(I_{2}+\gamma(y_2,y_3)\big)}\!-\!1.
\end{eqnarray*}
Similarly, define
\begin{equation}
\label{Delta(m)}
\Delta(y_1,y_2,\dots,y_m):=
\frac{{\rm det}
\big(I_m+\gamma(y_1,\dots, y_m)\big)}
{{\rm det}\big(I_{m-1}+\gamma(y_2,\dots,y_m)\big)}-1. 
\end{equation}
\begin{lem}[\cite{Kos-hpl-arx23}]
\label{l.m}
For the general $m$ we have
\begin{equation}
\label{m}
(C^{-1}(\lambda)a,a)=\Delta(y_1,y_2,\dots,y_m),
\end{equation}
where the $y_r$ are defined as follows:
\begin{equation}
\label{f_0,...f_m} 
y_r=
\left(\frac{a_{rk}}{\sqrt{\lambda_k}}\right)_{k=1}^n\in \mathbb R^n,\quad 1\leq r\leq m.
\end{equation} 
\end{lem}
%
\noindent{\it Acknowledgement.} 
The author is very grateful to Prof. K.-H. Neeb, Prof. M.~Smirnov and 
Dr Moree
for their personal efforts  to make academic  stays possible at their respective institutes. A.~Kosyak  visited: MPIM from March to April 2022 and from January to April 2023,
University of Augsburg from June to July 2022, and  University of Erlangen-Nuremberg 
from August to December 2022, all during the Russian invasion in Ukraine.
Also, Prof. R.~Kashaev  kindly invited A.~Kosyak to Geneva.

Further, he 
would like  to pay his respect to Prof. P. Teichner at MPIM, for his
immediate efforts started to help mathematicians
in Ukraine  after the Russian invasion.
Also Danil Radchenko 
helped the author 
to understand better the intial problem.

Since the spring of 2023 A.~Kosyak is an Arnold Fellow at the
London Institute for Mathematical Sciences, and he would like to express
his gratitude to  Mrs Myers Cornaby   to Miss Ker Mercer and to Miss Madeleine Hall 
and especially to the Director of LIMS Dr T.~Fink and Prof. Y.-H.~He.

\end{document}